\documentclass[12pt,draft]{amsart}
\usepackage{amsmath,amssymb,amsthm,bm}
\usepackage{graphicx,enumerate}
\usepackage{layout}

\theoremstyle{plain}
\newtheorem{theorem}{Theorem}[section]
\newtheorem{lemma}[theorem]{Lemma}

\newtheorem{corollary}[theorem]{Corollary}

\theoremstyle{definition}

\newtheorem{example}[theorem]{Example}
\newtheorem{remark}[theorem]{Remark}
\newtheorem*{acknowledgments}{Acknowledgments}
\theoremstyle{remark}

\newcommand{\bZ}{\mathbb{Z}}\newcommand{\bN}{\mathbb{N}}
\newcommand{\bQ}{\mathbb{Q}}\newcommand{\bR}{\mathbb{R}}
\newcommand{\cM}{\mathcal{M}}\newcommand{\cH}{\mathcal{H}}
\newcommand{\cC}{\mathcal{C}}\newcommand{\cR}{\mathcal{R}}
\newcommand{\cE}{\mathcal{E}}\newcommand{\cL}{\mathcal{L}^\natural}
\newcommand{\cO}{\mathcal{O}}
\newcommand{\iPi}{{\it\Pi}}
\newcommand{\bep}{\bm\varepsilon}\newcommand{\ep}{\varepsilon}

\newcommand{\te}{\theta}\newcommand{\bte}{\bm\theta}
\newcommand{\de}{\delta}\newcommand{\bde}{\bm\delta} 
\newcommand{\et}{\eta}
\newcommand{\bz}{\bm z}\newcommand{\bu}{\bm u}

\title[Cubic Thue equations and the simplest cubic fields]
{On correspondence between solutions of a family of cubic Thue 
equations and isomorphism classes of the simplest cubic fields}
\author{Akinari Hoshi}

\begin{document}
\maketitle
\begin{abstract}
Let $m\geq -1$ be an integer. 
We give a correspondence between integer solutions to the parametric family of 
cubic Thue equations
\[
X^3-mX^2Y-(m+3)XY^2-Y^3=\lambda
\]
where $\lambda>0$ is a divisor of $m^2+3m+9$ and isomorphism classes of 
the simplest cubic fields. 
By the correspondence and R. Okazaki's result, we determine the exactly 66 non-trivial 
solutions to the Thue equations for positive divisors $\lambda$ of $m^2+3m+9$. 
As a consequence, we obtain another proof of Okazaki's theorem which asserts that 
the simplest cubic fields are non-isomorphic to each other except for 
$m=-1,0,1,2,3,5,12,54,66,1259,2389$. 
\end{abstract}

\section{Introduction}\label{seIntro}

We consider the following parametric family of cubic Thue equations
\begin{align}
F_m(X,Y):=X^3-mX^2Y-(m+3)XY^2-Y^3=\lambda\label{eqFm}
\end{align}
for integer $m\in\bZ$ and non-zero integer $\lambda\in\bZ\setminus\{0\}$. 
We may assume that $m\geq -1$ and $\lambda>0$ because if $F_m(x,y)=\lambda$ then 
$F_m(-x,-y)=-\lambda$ and $F_{-m-3}(-y,-x)=\lambda$. 
For cubic integer $\lambda=c^3\in\bZ$, the equation (\ref{eqFm}) has three solutions 
\[
(c,0),\quad (0,-c),\quad (-c,c).
\]
In this paper, we call solutions $(x,y)\in\bZ^2$ to (\ref{eqFm}) 
with $xy(x+y)=0$ the \textit{trivial} solutions. 

In the case where $\lambda=\pm 1$, Thomas \cite{Tho90} showed that for 
$m':=m+1\geq 1.365\cdot 10^7$, the equation $F_m(x,y)=\pm 1$ has only the trivial 
solutions $(x,y)=(\pm 1,0)$, $(0,\mp 1)$, $(\mp 1,\pm 1)$ and for 
$0\leq m'=m+1\leq 10^3$, non-trivial solutions exist only for $m=-1, 0, 2$ (i.e. $m'=0,1,3$). 
Mignotte \cite{Mig93} solved completely the equations $F_m(X,Y)=\pm 1$ with the aid of 
a result of \cite{Tho90}. 
He proved that for $m\geq -1$, non-trivial solutions exist only for $m=-1,0,2$. 

For $m\geq -1$ and $\lambda=1$, all solutions to $F_m(x,y)=1$ are given by three trivial 
solutions $(x,y)=(1,0)$, $(0,-1)$, $(-1,1)$ for an arbitrary $m$ and additionally 
\begin{align*}
(x,y)&=(-1,-1),(-1,2),(2,-1),(5,4),(4,-9),(-9,5)\hspace*{-13mm}& &\mathrm{for}\quad m=-1,\\
(x,y)&=(2,1),(1,-3),(-3,2)& &\mathrm{for}\quad m=0,\\
(x,y)&=(-7,-2),(-2,9),(9,-7)& &\mathrm{for}\quad m=2
\end{align*}
(cf. also \cite[page 54]{Gaa02}). 
Note that if $(x,y)\in\bZ^2$ is a solution to (\ref{eqFm}) then $(y,-x-y)$ and $(-x-y,x)$ 
are also solutions to (\ref{eqFm}) because $F_m(X,Y)$ is invariant under the action of 
the cyclic group $C_3=\langle\sigma\rangle$ of order three 
where $\sigma\,:\, X\longmapsto Y\longmapsto -X-Y$. 

Mignotte-Peth\"o-Lemmermeyer \cite{MPL96} studied the equation (\ref{eqFm}) 
for general $\lambda\in\bZ$ and gave a complete solution to Thue inequality 
$|F_m(X,Y)|\leq 2m+3=2m'+1$ with a result of \cite{LP95}. 
For $m\geq -1$ and $1<\lambda\leq 2m+3$, all solutions to (\ref{eqFm}) are given by 
trivial solutions for $\lambda=c^3$ and 
\begin{align*}
(x,y)\in\{(-1,-1),(-1,2),(2,-1),(-m-1,-1),(-1,m+2),(m+2,-m-1)\}
\end{align*}
for $\lambda=2m+3$, except for $m=1$ in which case (\ref{eqFm}) has the extra solutions: 
\begin{align*}
(x,y)\in\{(3,1),(1,-4),(-4,3),(8,3),(3,-11),(-11,8)\}
\end{align*}
for $\lambda=5$
\footnote{In \cite[Theorem 3]{MPL96}, there is a misprint: it should be added $(-11,8)$.}.
Lettl-Peth\"o-Voutier \cite{LPV99} and Xia-Chen-Zhang \cite{XCZ06} also investigated Thue 
inequality $|F_m(X,Y)|\leq \lambda(m)$ where $\lambda\,:\,\bZ\rightarrow \bN$ under 
some condition. 
Wakabayashi \cite{Wak07b} studied Thue inequality $|F_{l,m}(X,Y)|\leq\lambda$ 
with two parameters $l,m$ and $F_{1,m}=F_m$. 
We will use a result of \cite{LPV99} in Section \ref{seThue}.

Let $\bQ(t)$ be the rational function field over $\bQ$ with variable $t$. 
We take a polynomial 
\[
f_t^{C_3}(X):=F_t(X,1)=X^3-tX^2-(t+3)X-1\in \bQ(t)[X]
\]
with discriminant $(t^2+3t+9)^2$. 
The polynomial $f_t^{C_3}(X)$ is generic for $C_3$ over $\bQ$ in the sense of \cite{JLY02}. 
Namely, the Galois group of $f_t^{C_3}(X)$ over $\bQ(t)$ is isomorphic to $C_3$ and every 
Galois extension $L/M$, $M\supset\bQ$, can be obtained as $L=\mathrm{Spl}_M f_a^{C_3}(X)$, 
the splitting field of $f_a^{C_3}(X)$ over $M$, for some $a\in M$. 
For $m\in\bZ$, the polynomials $f_m^{C_3}(X)$ are irreducible over $\bQ$ and 
the splitting fields 
\[
L_m:=\mathrm{Spl}_\mathbb{Q} f_m^{C_3}(X)
\]
of $f_m^{C_3}(X)$ over $\mathbb{Q}$ are called the simplest cubic fields (cf. \cite{Sha74}). 

Because $F_m(X,1)=F_{-m-3}(-1,-X)$, we see if $z$ is a root of $f_m^{C_3}(X)$ then 
$1/z$ becomes a root of $f_{-m-3}^{C_3}(X)$. 
Hence $L_m=L_{-m-3}$, for any $m\in\bZ$. 

The aim of this paper is to give the following theorem which asserts a correspondence 
between certain non-trivial solutions to (\ref{eqFm}) and isomorphic simplest cubic fields. 

\begin{theorem}\label{thmain}
Let $m$ be an integer and $L_m$ the splitting field of $f_m^{C_3}(X)
=X^3-mX^2-(m+3)X-1$ over $\bQ$. 
There exists an integer $n\in\bZ\setminus\{m,-m-3\}$ such that 
$L_n=L_m$ if and only if there exists a solution $(x,y)\in\bZ^2$ with 
$xy(x+y)\neq 0$ to the family of cubic Thue equations 
\begin{align}
F_m(x,y)=x^3-mx^2y-(m+3)xy^2-y^3=\lambda\label{eqFth}
\end{align}
where $\lambda>0$ is a divisor of $m^2+3m+9$. 
Moreover, an integer $n$ and solutions $(x,y)\in\bZ^2$ to $(\ref{eqFth})$ can be chosen 
to satisfy 
\begin{align}
N=m+\frac{(m^2+3m+9)xy(x+y)}{F_m(x,y)}\label{eqth2}
\end{align}
where either $N=n$ or $N=-n-3$, and this occurs for only one of $n$ and $-n-3$. 
\end{theorem}
\begin{corollary}\label{corC3}
Let $\mathcal{N}$ be the number of primitive solutions $(x,y)\in\bZ^2$, 
i.e. $\mathrm{gcd}(x,y)=1$, with $xy(x+y)\neq0$ to $(\ref{eqFth})$  for a fixed $m\in\bZ$. 
For $m\in\bZ$, we have 
\begin{align*}
\#\bigl\{n\in\bZ\setminus\{m,-m-3\}\ \big|\ L_n=L_m,\, n\geq -1\bigr\}=\frac{\mathcal{N}}{3}.
\end{align*}
\end{corollary}

Note that (i) the discriminant of $f_m^{C_3}(X)$ with respect to $X$ 
equals $(m^2+3m+9)^2$, and 
(ii) if $(x,y)\in\bZ^2$ is a solution to $(\ref{eqFth})$  then $(y,-x-y)$, $(-x-y,x)$ 
are also solutions to $(\ref{eqFth})$. 

Ennola \cite{Enn91} verified that for integers $-1\leq m<n\leq 10^4$, 
the overlap $L_m=L_n$ of the splitting fields occurs 
if and only if $(m,n)\in \{(-1,5)$, $(-1,12)$, $(-1,1259)$, $(5,12)$, $(5,1259)$, 
$(12,1259)\}$ $\cup$ $\{(0,3)$, $(0,54)$, $(3,54)\}$ $\cup$ $\{(1,66)\}$ $\cup$ $\{(2,2389)\}$. 
Hoshi-Miyake \cite{HM09a} checked with the aid of computer that this claim is also valid for 
$-1\leq m<n\leq 10^5$ (see \cite[Example 5.3]{HM09a}). 

Okazaki \cite{Oka02} investigated Thue equations $F(X,Y)=1$ for irreducible cubic forms $F$ 
with positive discriminant $D(F)>0$ and established a very strong result on gaps between 
solutions (cf. Theorem \ref{thOkaG}). 
Using results in \cite{Oka02}, he showed the following theorems: 
\begin{theorem}[Okazaki]\label{thOkab}
For integers $-1\leq m<n$, if $L_m=L_n$ then $m\leq 35731$. 
\end{theorem}
\begin{theorem}[Okazaki]\label{thOka}
For integers $-1\leq m<n$, if $L_m=L_n$ then 
$m,n\in\{-1$, $0$, $1$, $2$, $3$, $5$, $12$, $54$, $66$, $1259$, $2389\}$. In particular, 
\[
L_{-1}=L_5=L_{12}=L_{1259},\quad L_0=L_3=L_{54},\quad L_1=L_{66},\quad L_2=L_{2389}. 
\]
\end{theorem}
Although Theorems \ref{thOkab} and \ref{thOka} seem to be unpublished, 
a brief sketch of the proof is available 
at \verb+http://www1.doshisha.ac.jp/~rokazaki/papers.html+ 
as a presentation sheet \cite{Oka} (cf. also \cite{Wak07a}). 
In Section \ref{seOka}, we will recall Okazaki's results in \cite{Oka02} 
and give a proof of Theorem \ref{thOkab}. 
Theorem \ref{thOkab} can be shown by the same method as in \cite{Oka02}. 
In Section \ref{seThue}, we determine all solutions to (\ref{eqFth}) by using 
Theorems \ref{thmain} and \ref{thOkab}, Bilu-Hanrot \cite{BH96} and 
Lettl-Peth\"o-Voutier \cite{LPV99}. 
Then we get another proof of Theorem \ref{thOka} as a consequence of 
Theorems \ref{thmain} and \ref{thallsol}. 
\begin{theorem}\label{thallsol}
For $m\geq -1$, there exist exactly $66$ solutions $(x,y)\in\bZ^2$ with $xy(x+y)\neq 0$ 
to the family of Thue equations $F_m(x,y)=\lambda$, where $\lambda>0$ is a divisor of $m^2+3m+9$. 
These $66$ solutions are given in Table 1. 
\end{theorem}
\begin{center}
{\rm Table} $1$\vspace*{4mm}\\
{\scriptsize
\begin{tabular}{|c|c|c|c|c|c|c|l|}\hline
$m$ & $N$ & $-N-3$ & $2m+3$ & $\lambda$ & $m^2+3m+9$ & $xy(x+y)$ & $(x,y)$\\ \hline 
$-1$ & $-15$ & $12$ & $1$ & $1$ & $7$ & $-2$ & $(-1,-1)$, $(-1,2)$, $(2,-1)$\\ \hline
$-1$ & $1259$ & $-1262$ & $1$ & $1$ & $7$ & $180$ & $(5,4)$, $(4,-9)$, $(-9,5)$\\ \hline
$-1$ & $5$ & $-8$ & $1$ & $7$  & $7$ & $6$ & $(2,1)$, $(1,-3)$, $(-3,2)$\\ \hline
$0$ & $54$ & $-57$ & $3$ & $1$ & $9$ & $6$ & $(2,1)$, $(1,-3)$, $(-3,2)$\\ \hline 
$0$ & $-6$ & $3$ & $3$ & $3$ & $9$ & $-2$ & $(-1,-1)$, $(-1,2)$, $(2,-1)$\\ \hline 
$1$ & $-69$ & $66$ & $5$ & $13$ & $13$ & $-70$ & $(-5,-2)$, $(-2,7)$, $(7,-5)$\\ \hline 
$2$ & $-2392$ & $2389$ & $7$ & $1$ & $19$ & $-126$ & $(-7,-2)$, $(-2,9)$, $(9,-7)$\\ \hline 
$3$ & $-3$ & $0$ & $9$ & $9$ & $3^3$ & $-2$ & $(-1,-1)$, $(-1,2)$, $(2,-1)$\\ \hline 
$3$ & $-57$ & $54$ & $9$ & $9$ & $3^3$ & $-20$ & $(-4,-1)$, $(-1,5)$, $(5,-4)$\\ \hline 
$5$ & $-1$ & $-2$ & $13$ & $7^2$ & $7^2$ & $-2$ & $(-1,-2)$, $(-2,3)$, $(3,-1)$\\ \hline 
$5$ & $-15$ & $12$ & $13$ & $7^2$ & $7^2$ & $-20$ & $(-4,-1)$, $(-1,5)$, $(5,-4)$\\ \hline 
$5$ & $1259$ & $-1262$ & $13$ & $7^2$ & $7^2$ & $1254$ & $(19,3)$, $(3,-22)$, $(-22,19)$\\ \hline 
$12$ & $-2$ & $-1$ & $27$ & $3^3$ & $3^3\cdot 7$ & $-2$ & $(-1,-1)$, $(-1,2)$, $(2,-1)$\\ \hline 
$12$ & $-1262$ & $1259$ & $27$ & $3^3$ & $3^3\cdot 7$ & $-182$ 
& $(-13,-1)$, $(-1,14)$, $(14,-13)$\\ \hline 
$12$ & $-8$ & $5$ & $27$ & $189=3^3\cdot 7$ & $3^3\cdot 7$ & $-20$ 
& $(-4,-1)$, $(-1,5)$, $(5,-4)$\\ \hline 
$54$ & $0$ & $-3$ & $111$ & $343=7^3$ & $3^2\cdot 7^3$ & $-6$ 
& $(-1,-2)$, $(-2,3)$, $(3,-1)$\\ \hline 
$54$ & $-6$ & $3$ & $111$ & $1029=3\cdot 7^3$ & $3^2\cdot 7^3$ & $-20$ 
& $(-4,-1)$,$(-1,5)$, $(5,-4)$\\ \hline 
$66$ & $-4$ & $1$ & $135$ & $4563=3^3\cdot 13^2$ & $3^3\cdot 13^2$ & $-70$ 
& $(-5,-2)$, $(-2,7)$, $(7,-5)$\\ \hline 
$1259$ & $-1$ & $-2$ & $2521$ & $226981=61^3$ & $7\cdot 61^3$ & $-180$ 
& $(-4,-5)$, $(-5,9)$, $(9,-4)$\\ \hline
$1259$ & $-15$ & $12$ & $2521$ & $226981=61^3$ & $7\cdot 61^3$ & $-182$ 
& $(-13,-1)$, $(-1,14)$, $(14,-13)$\\ \hline
$1259$ & $5$ & $-8$ & $2521$ & $1588867=7\cdot 61^3$ & $7\cdot 61^3$ & $-1254$ 
& $(-3,-19)$, $(-19,22)$, $(22,-3)$\\ \hline 
$2389$ & $-5$ & $2$ & $4781$ & $300763=67^3$ & $19\cdot 67^3$ & $-126$ 
& $(-7,-2)$, $(-2,9)$, $(9,-7)$\\ \hline 
\end{tabular}
}
\vspace*{1mm}
\end{center}

In Theorem \ref{thallsol}, only primitive solution $(x,y)\in\bZ^2$, i.e. $\gcd(x,y)=1$, 
to $F_m(x,y)=\lambda$ exists (see Table $1$), since $(m^2+3m+9)/\lambda$ are cubefree. 
\begin{corollary}\label{corsol}
For $m\geq -1$, the only solutions $(x,y)\in\bZ^2$ to the parametric family of 
Thue equations $F_m(x,y)=m^2+3m+9$ are given by
\begin{align*}
(x,y)=&\ (2,1),(1,-3),(-3,2)& &\hspace*{-15mm}\mathrm{for}\quad m=-1,\\
(x,y)=&\ (-5,-2),(-2,7),(7,-5)& &\hspace*{-15mm}\mathrm{for}\quad m=1,\\
(x,y)=&\ (3,0),(0,-3),(-3,3)& &\hspace*{-15mm}\mathrm{for}\quad m=3,\\
(x,y)=&\ (-1,-2),(-2,3),(3,-1),(-4,-1),& &\\
&\ (-1,5),(5,-4),(19,3),(3,-22),(-22,19)& &\hspace*{-15mm}\mathrm{for}\quad m=5,\\
(x,y)=&\ (-4,-1),(-1,5),(5,-4)& &\hspace*{-15mm}\mathrm{for}\quad m=12,\\
(x,y)=&\ (-5,-2),(-2,7),(7,-5)& &\hspace*{-15mm}\mathrm{for}\quad m=66,\\
(x,y)=&\ (-3,-19),(-19,22),(22,-3)& &\hspace*{-15mm}\mathrm{for}\quad m=1259.
\end{align*}
\end{corollary}
\begin{remark}
Theorem \ref{thOka} may be obtained as a consequence of 
Theorems \ref{thmain} and \ref{thallsol}. 
Conversely if we assume Theorem \ref{thOka} then 
by applying Theorem \ref{thmain} to fixed $m$ and $n$, 
we may get all solutions to $(\ref{eqFth})$  as on Table $1$ without the aid of computer. 
\end{remark}
\section{Field isomorphism problem}\label{sePre}

In order to prove Theorem \ref{thmain}, 
we recall a result of \cite{Mor94}, \cite{Cha96} and \cite{HM09a} 
for the field isomorphism problem of the simplest cubic polynomial $f_t^{C_3}(X)$, i.e. 
for a fixed $m,n\in K$ whether the splitting fields of $f_m^{C_3}(X)$ and of $f_n^{C_3}(X)$ 
over $K$ coincide. 

Let $K$ be a field of char $K\neq 2$ and $K(t)$ the rational function field over $K$ 
with variable $t$. 
We take the simplest cubic polynomial 
\[
f_t^{C_3}(X)=X^3-tX^2-(t+3)X-1\in K(t)[X].
\]
Note that the splitting fields of $f_m^{C_3}(X)$ and of $f_{-m-3}^{C_3}(X)$ over $K$ 
coincide for any $m\in K$. 

In \cite{HM09a}, by using formal Tschirnhausen transformation, we showed the following theorem 
which gives an answer to the field isomorphism problem of $f_t^{C_3}(X)$. 
This theorem was given by Morton \cite{Mor94} essentially (see also \cite{Cha96} 
and \cite{HM09a} with replacing $z$ by $1/z$). 

\begin{theorem}[{\cite[Theorem 3]{Mor94}, \cite[Corollary 1]{Cha96}, 
\cite[Theorem 5.4]{HM09a}}]\label{thC3}
Let $K$ be a field of char $K\neq 2$. 
For $m\in K$ and $n\in K\setminus\{m,-m-3\}$ with $(m^2+3m+9)(n^2+3n+9)\neq 0$, 
the following three conditions are equivalent{\rm :}\\
{\rm (i)} the splitting fields of $f_m^{C_3}(X)$ and of $f_n^{C_3}(X)$ over $K$ coincide{\rm ;}\\
{\rm (ii)} the polynomial $f_{M_i}^{C_3}(X)$ splits completely into three linear factors 
over $K$ for $i=1$ or $i=2$ 
where
\[
M_1=-\frac{mn+3m+9}{m-n}\quad\textrm{or}\quad M_2=\frac{mn-9}{m+n+3};
\]
{\rm (iii)} there exists $z\in K$ such that 
\[
N\,=\,m+\frac{(m^2+3m+9)z(z+1)}{f_m^{C_3}(z)}
\]
where $N=n$ or $N=-n-3$.

Moreover, the conditions {\rm (i)} and {\rm (iii)} are equivalent also for $n=m$ and $n=-m-3$, 
and if $\mathrm{Gal}_K f_m^{C_3}(X)\cong C_3$ 
$($resp. $\mathrm{Gal}_K f_m^{C_3}(X)\cong \{1\})$ then 
{\rm (ii)} occurs for only one of $M_1$ and $M_2$ $($resp. for both of $M_1$ and $M_2)$ and 
{\rm (iii)} occurs for only one of $n$ and $-n-3$ $($resp. for both of $n$ and $-n-3)$. 
\end{theorem}
\begin{proof}
Here we give a sketch of proof. 
For simplicity, we assume that 
$\mathrm{Gal}_K f_m^{C_3}(X)\cong \mathrm{Gal}_K f_n^{C_3}(X)\cong C_3$ 
(see \cite{HM09a} and \cite{HM09b} for general case). 
We first note that {\rm (iii)} is just a restatement of {\rm (ii)}. 
Take a root $\te$ of $f_m^{C_3}(X)$ and a root $\et$ of $f_n^{C_3}(X)$ in $\overline{K}$, 
then we have $m=(\te^3-3\te-1)/(\te(\te+1))$, $n=(\et^3-3\et-1)/(\et(\et+1))$, 
$\mathrm{Gal}_K K(\te)=\langle\tau_1\rangle$, $\mathrm{Gal}_K K(\et)=\langle\tau_2\rangle$ and 
$K(\te,\et)^{\langle\tau_1\rangle\times\langle\tau_2\rangle}=K(m,n)=K$
where 
\begin{align*}
\tau_1\,:\, \te\mapsto \frac{1}{-\te-1}\mapsto \frac{-\te-1}{\te}\mapsto \te,\qquad 
\tau_2\,:\, \et\mapsto \frac{1}{-\et-1}\mapsto \frac{-\et-1}{\et}\mapsto \et. 
\end{align*}
Now we put 
\[
\Theta_1:=-\frac{\te\et+\te+1}{\te-\et},\qquad \Theta_2:=\frac{\te\et-1}{\te+\et+1}.
\] 
Then we get $K(\te,\et)^{\langle\tau_1\tau_2\rangle}=K(\Theta_1)$ and 
$K(\te,\et)^{\langle\tau_1\tau_2^2\rangle}=K(\Theta_2)$. 
The orbit of $\Theta_i$, $(i=1,2)$, under the action of 
$\mathrm{Gal}_K K(\te,\et)=\langle\tau_1\rangle\times\langle\tau_2\rangle$ 
is given by the same as the $\langle\tau_1\rangle$-orbit of $\te$: 
\[
\mathrm{Orb}_{\langle\tau_1\rangle\times\langle\tau_2\rangle}(\Theta_i)
=\Big\{\Theta_i,\ \frac{1}{-\Theta_i-1},\ \frac{-\Theta_i-1}{\Theta_i}\Big\}. 
\]
Therefore the minimal polynomial of $\Theta_i$, $(i=1,2)$, over $K$ is also given as 
``the simplest cubic polynomial'' $f_M^{C_3}(X)$ for some $M\in K(m,n)=K$. 
Indeed we may evaluate such $M$'s as in the condition {\rm (ii)} respectively. 

It follows from the theory of resolvent polynomials in Galois theory 
that the condition {\rm (i)} holds if and only if either $\Theta_1\in K$ or $\Theta_2\in K$ 
(cf. for example \cite[Chapter 6]{Ade01}, \cite{HM09b}, \cite{HM09c}, \cite{HM10}). 
\end{proof}
We consider the case where $K=\bQ$. 
Let $m\in\bZ$ and $L_m=\mathrm{Spl}_\bQ f_m^{C_3}(X)$. 
The discriminant of $f_m^{C_3}(X)$ is $(m^2+3m+9)^2$. 
We note that $3\,|\,m^2+3m+9$ if and only if $3\,|\,m$. 
We also see $3^2\,||\,m^2+3m+9$ if $m\equiv0,6\pmod{9}$ and 
$3^3\,||\,m^2+3m+9$ if $m\equiv3\pmod{9}$ where $3^r\,||\,a$ 
means $3^r\mid a$ and $3^{r+1}\!{\not{\mid}}\ a$. 
The conductor $\mathfrak{f}_m$ and the discriminant $d_{L_m}=\mathfrak{f}_m^2$ 
of the field $L_m$ are given as follows 
(cf. e.g. \cite{Gra73}, \cite[Proposition 3.10]{Gra86}, 
\cite[Proposition 1]{Was87}, \cite[Theorem 3.6]{Kom04}, 
\cite[Lemma 1.4]{Kom07}, \cite[Theorem 10]{Hab10}). 
\begin{lemma}\label{lemcond}
Let $p_1,\ldots,p_s$ be primes, different from $3$, dividing $m^2+3m+9$ 
with an exponent not congruent to $0$ modulo $3$. 
Then the conductor $\mathfrak{f}_m$ of $L_m$ is given by 
\begin{align*}
\mathfrak{f}_m=
\begin{cases}
p_1\cdots p_s\hspace*{6.7mm}\mathrm{if}\ m\not\equiv0\!\!\!\pmod{3}\ 
\mathrm{or}\ m\equiv12\!\!\!\pmod{27},\\
3^2 p_1\cdots p_s\ \ \mathrm{if}\ m\equiv0\!\!\!\pmod{3}\ 
\mathrm{and}\ m\not\equiv12\!\!\!\pmod{27}.
\end{cases}
\end{align*}
\end{lemma}
In Section \ref{seThue}, we will use Lemma \ref{lemcond} to prove Corollary \ref{corsol}. 
Although there exist $2^{s-1}$ (resp. $2^s$) cyclic cubic fields 
with conductor $\mathfrak{f}_m$, by Lemma \ref{lemcond}, we may confirm that 
$L_m\neq L_n$ if $\mathfrak{f}_m\neq \mathfrak{f}_n$. 
When $\mathfrak{f}_m=\mathfrak{f}_n$, we need Theorem \ref{thC3} to decide 
whether $L_m=L_n$ or not. 

\begin{example} 
We give numerical examples which satisfy $\mathfrak{f}_m=\mathfrak{f}_n$ and $L_m\neq L_n$ 
with $-1\leq n<m\leq 1000$ as on Table $2$. 

\begin{center}
{\rm Table} $2$\vspace*{4mm}\\
\begin{tabular}{|c|c|c|c|c|}\hline
$m$ & $m^2+3m+9$ & $n$ & $n^2+3n+9$ & $\mathfrak{f}_m=\mathfrak{f}_n$ \\ \hline 
$13$ & $7\cdot 31$ & $201$ & $3^3\cdot 7^2\cdot 31$ & $217=7\cdot 31$\\\hline
$27$ & $3^2\cdot 7\cdot 13$ & $48$ & $3^3\cdot 7\cdot 13$ & $819=3^2\cdot 7\cdot 13$\\\hline
$33$ & $3^2\cdot 7\cdot 19$ & $90$ & $3^2\cdot 7^2\cdot 19$ & $1197=3^2\cdot 7\cdot 19$\\\hline
$51$ & $3^2\cdot 307$ & $972$ & $3^2\cdot 7^3\cdot 307$ & $2763=3^2\cdot 307$\\\hline
$53$ & $13\cdot 229$ & $282$ & $3^3\cdot 13\cdot 229$ & $2977=13\cdot 229$\\\hline
$79$ & $13\cdot 499$ & $417$ & $3^3\cdot 13\cdot 499$ & $6487=13\cdot 499$\\\hline
$105$ & $3^2\cdot 13\cdot 97$ & $183$ & $3^3\cdot13\cdot 97$ & 
$11349=3^2\cdot 13\cdot 97$\\\hline
$261$ & $3^2\cdot 13\cdot 19\cdot 31$ & $945$ & $3^2\cdot 13^2\cdot 19\cdot 31$ & 
$68913=3^2\cdot13\cdot19\cdot31$\\\hline
$516$& $3^3\cdot7\cdot13\cdot109$ & $789$ & $3^2\cdot7^2\cdot13\cdot109$ & 
$89271=3^2\cdot7\cdot13\cdot109$\\\hline
\end{tabular}
\end{center}
\vspace*{2mm}
\end{example}

%%%%%%%%%%%%%%%%%%%%%%%%%%%%%%%%%%%%%%%%%%%%%%%%%%%%%%%%%%%%%%%%%%%%%%%%%%%%%%%%%%
\begin{proof}[Proof of Theorem \ref{thmain}]
We use Theorem \ref{thC3} in the case where $K=\bQ$ and $m,n\in\bZ$.
Assume that there exists an integer $n\in\bZ\setminus\{m,-m-3\}$ such that $L_m=L_n$. 
By Theorem \ref{thC3}, there exist $x,y\in\bZ$ with $z=x/y$ and $xy(x+y)\neq 0$ such that
\begin{align*}
{N}=m+\frac{(m^2+3m+9)xy(x+y)}{F_m(x,y)}
\end{align*}
where either $N=n$ or $N=-n-3$. 
In particular, we have 
\begin{align}
\frac{(m^2+3m+9)xy(x+y)}{F_m(x,y)}\in\bZ.\label{eqInt}
\end{align}
We will show that if we take $(x,y)\in\bZ^2$ with $z=x/y$ as $\gcd(x,y)=1$ then 
$\lambda:=F_m(x,y)$ divides $m^2+3m+9$. 
Put $h(z):=(m^2+3m+9)z(z+1)$. 
We take the resultant $R:=\mathrm{Res}_z(h(z),f_m^{C_3}(z))$ of $h(z)$ and $f_m^{C_3}(z)$ 
with respect to $z$. 
By definition, 
\begin{align}
R={\footnotesize 
\left| \begin{array}{ccccc}
m^2+3m+9 & m^2+3m+9 & 0 & 0 & 0\\
0 & m^2+3m+9 & m^2+3m+9  & 0 & 0\\
0 & 0 & m^2+3m+9 & m^2+3m+9  & 0 \\
 1 & -m & -m-3 & -1 & 0 \\
 0 & 1 & -m & -m-3 & -1
\end{array} \right|}=-(m^2+3m+9)^3\label{R1}.
\end{align}
By matrix computations, we see easily that $R$ is also given as of the form 
\begin{align}
R&={\small 
\left| \begin{array}{cccccc}
m^2+3m+9 & m^2+3m+9 & 0 & 0 & h(z)z^2\\
0 & m^2+3m+9 & m^2+3m+9  & 0 & h(z)z\\
0 & 0 & m^2+3m+9 & m^2+3m+9  & h(z) \\
1 & -m & -m-3 & -1 & f_m^{C_3}(z)z \\
0 &1 & -m & -m-3 & f_m^{C_3}(z) \\
\end{array} \right|}\nonumber\\
&=-(m^2+3m+9)^2\Bigl(h(z)p(z)+f_m^{C_3}(z)q(z)\Bigr)\label{R11}
\end{align}
where
\begin{align*}
p(z)=2z^2-2mz-z-m-5,\quad q(z)=-(m^2+3m+9)(2z+1).
\end{align*}
Then, by (\ref{R1}) and (\ref{R11}), we have 
\begin{align*}
h(z)p(z)+f_m^{C_3}(z)q(z)=m^2+3m+9.
\end{align*}
Put 
\begin{align*}
H(x,y)&:=(m^2+3m+9)xy(x+y),\\
P(x,y)&:=2x^2-2mxy-xy-my^2-5y^2,\\
Q(x,y)&:=-(m^2+3m+9)y(2x+y). 
\end{align*}
Then it follows from $z=x/y$ that 
\begin{align*}
H(x,y)P(x,y)+F_m(x,y)Q(x,y)=(m^2+3m+9)y^5.
\end{align*}
Because the cubic forms $F_m(X,Y)$ and $H(X,Y)$ are invariants under the action of 
$\sigma^2\,:\, Y\longmapsto X\longmapsto -X-Y$, we also get 
\begin{align*}
H(x,y)P(-x-y,x)+F_m(x,y)Q(-x-y,x)=(m^2+3m+9)x^5.
\end{align*}
Hence by (\ref{eqInt}) we have
\begin{align*}
&\frac{H(x,y)P(x,y)}{F_m(x,y)}+Q(x,y)=\frac{(m^2+3m+9)y^5}{F_m(x,y)}\in\bZ,\\
&\frac{H(x,y)P(-x-y,x)}{F_m(x,y)}+Q(-x-y,x)=\frac{(m^2+3m+9)x^5}{F_m(x,y)}\in\bZ.
\end{align*}
By the assumption $\mathrm{gcd}(x,y)=1$, we conclude that $\lambda=F_m(x,y)$ divides $m^2+3m+9$. 

Conversely if there exists $(x,y)\in\bZ^2$ such that 
$\lambda=F_m(x,y)>0$ divides $m^2+3m+9$ then we may choose $(x,y)$ with 
$\mathrm{gcd}(x,y)=1$ and also get $n\in\bZ$ which satisfies $L_n=L_m$ by
\begin{align}
n=m+\frac{(m^2+3m+9)xy(x+y)}{F_m(x,y)}.\label{eqmmm}
\end{align}
Note that $n=m$ if and only if $xy(x+y)=0$. 
By Theorem \ref{thC3}, there does not exist $(x,y)\in\bZ^2$ 
which satisfies (\ref{eqmmm}) for $n=-m-3$  because 
$\mathrm{Gal}_\bQ f_m^{C_3}(X)\cong C_3$. 
Hence if $(x,y)\in\bZ^2$ with $xy(x+y)\neq 0$ satisfies (\ref{eqmmm}) then 
$n\in\bZ\setminus\{m,-m-3\}$. 
\end{proof}

\begin{proof}[Proof of Corollary \ref{corC3}]
By Theorem \ref{thC3}, if $(x,y)\in\bZ^2$ with $\mathrm{gcd}(x,y)=1$ satisfies (\ref{eqmmm}) 
then all primitive solutions to $F_m(x,y)=\lambda>0$ which satisfy (\ref{eqmmm}) for 
the same $n$ are given by $(x,y)$, $(y,-x-y)$, $(-x-y,x)$. 
By Theorem \ref{thC3} again, there does not exist $(x,y)\in\bZ^2$ 
which gives $-n-3$ instead of $n$ via (\ref{eqmmm}) because 
$\mathrm{Gal}_\bQ f_m^{C_3}(X)\cong C_3$. 
\end{proof}

%
%
%%%%%%%%%%%%%%%%%%%%%%%%%%%%%%%%%%%%%%%%%%%%%%%%%%%%%%%%%%%%%%%%%%%%%%%%%%%%%%%%%%%%%%%%%
\section{Okazaki's results}\label{seOka}

In this section, we recall Okazaki's results in \cite{Oka02} 
(cf. also \cite{Wak07a}, \cite{Akh09}) 
and give a proof of Theorem \ref{thOkab}.

Let $F(X,Y)\in\bZ[X,Y]$ be an irreducible cubic form of positive discriminant 
$D=D(F)>0$. 
Let $\cR=\cR(F)$ be the set of integer solutions $(x,y)\in\bZ^2$ of
\begin{align*}
F(x,y)=1,
\end{align*}
i.e. the set of representations of $1$ by $F$ and $\mathrm{Aut}(F)$ be the group 
of automorphisms of $F$. 
In his paper \cite{Oka02}, 
Okazaki established a strong estimate for gaps between solutions of cubic Thue equations 
and gave an upper bound to $\#\cR$ under some conditions. 
In the case of $\mathrm{Aut}(F)\neq 1$, the estimation for $\#\cR$ becomes more efficient. 
Indeed, using his gap principle (see Theorem \ref{thOkaG} below), 
Okazaki showed the following theorem:  
\begin{theorem}[{\cite[Theorem 1.3]{Oka02}}]\label{thOkaM}
Assume $\mathrm{Aut}(F)\neq 1$. 
If $D(F)\geq 2.56\cdot 10^{18}$, then we have $\#\cR(F)\in\{0,3\}$.
\end{theorem}
This result gives a generalization of Thomas' result \cite{Tho90} which we mentioned 
in Section \ref{seIntro} (cf. also \cite{Mig98}). 
Indeed, for our case $F=F_m$, if $m\geq 4\cdot 10^4$ then 
$D(F_m)=(m^2+3m+9)^2> 2.56\cdot 10^{18}$ and hence $\#\cR(F_m)=3$, i.e. 
$F_m(X,Y)=1$ has only three trivial solutions (cf. Thomas' estimate $m+1\geq 1.365\cdot 10^7$). 

For our purpose, we restrict ourselves to consider only the case 
\begin{align*}
\mathrm{Aut}(F)\neq 1,\ \#\cR>0,\ \mathrm{and}\ 
f(X)=F(X,1)\ \mathrm{is\ monic\ and\ irreducible,}
\end{align*}
(see \cite{Oka02} for general case). 
Let $\te_1,\te_2,\te_3$ be roots of $f(X)$. 
Then $\bQ(\te_1,\te_2,\te_3)=\bQ(\te_1)$ is a totally real cyclic cubic field 
over $\bQ$ with $\mathrm{Gal}(\bQ(\te_1)/\bQ)=\langle\sigma\rangle$, 
\begin{align}
\sigma : \te_1\mapsto \te_2\mapsto \te_3\mapsto \te_1,\label{acts}
\end{align} and 
\begin{align*}
F(X,Y)=(X-\te_1 Y)(X-\te_2 Y)(X-\te_3 Y),\quad 
D=((\te_1-\te_2)(\te_1-\te_3)(\te_2-\te_3))^2.
\end{align*}
Let $\cO(F)=\bZ[\te_1]$ be an order in the ring of algebraic integers of $\bQ(\te_1)$ 
and $\cO(F)^\times$ the unit group in $\cO(F)$. 
For two vectors $\mathbf{1}={}^t(1,1,1)$ and $\bte={}^t(\te_1,\te_2,\te_3)$ in $\bR^3$, 
we take the exterior product 
\begin{align*}
\bde={}^t(\de_1,\de_2,\de_3):=\mathbf{1}\times\bte={}^t(\te_2-\te_3,\te_3-\te_1,\te_1-\te_2).
\end{align*}
Note that the norm $\mathrm{N}(\bde)=\de_1\de_2\de_3$ is given by $\mathrm{N}(\bde)=-\sqrt{D}$ 
and $\bde^{1-\sigma}\in (\cO(F)^\times)^3$ because $\mathrm{N}(\bde^{1-\sigma})=1$. 
The canonical lattice 
\[
\cL=\bde(\bZ\mathbf{1}+\bZ\bte)
\]
of $F$ is orthogonal to $\mathbf{1}$, where the product of vectors 
is the component-wise product. 
We consider the plane 
$
\iPi=\{{}^t(z_1,z_2,z_3)\in\bR^3 \mid z_1+z_2+z_3=0\}
$
and the curve $\cH$ on $\iPi$: 
\[
\cH : z_1+z_2+z_3=0,\quad z_1z_2z_3=\sqrt{D}.
\]
For $(x,y)\in\cR$, we see $x\mathbf{1}-y\bte\in (\cO(F)^\times)^3$ with 
$\mathrm{N}(x\mathbf{1}-y\bte)=1$. 
Then we get a bijection 
\[
\cR\ni(x,y)\longleftrightarrow \bz=\bde(-x\mathbf{1}+y\bte)\in\cL\cap\cH
\]
via $\mathrm{N}(\bz)=\mathrm{N}(\bde)\mathrm{N}(-x\mathbf{1}+y\bte)=(-\sqrt{D})(-1)=\sqrt{D}$. 
Then we have $\#\cR=\#(\cL\cap\cH)$. 

If $\bz\in\cL\cap\cH$ then 
$\bz^\sigma={}^t(z_1^\sigma,z_2^\sigma,z_3^\sigma)={}^t(z_2,z_3,z_1)$, 
$\bz^{\sigma^2}={}^t(z_3,z_1,z_2)\in\cL\cap\cH$. 
From this, we use the notation $\sigma$ both for $\sigma\in\mathrm{Gal}(\bQ(\theta_1)/\bQ)$ 
as in (\ref{acts}) and for the rotation 
$\sigma : {}^t(z_1,z_2,z_3)\mapsto {}^t(z_2,z_3,z_3)\in\bR^3$ 
(cf. \cite[Theorem 4.2]{Oka02}). 
In particular, $3\mid\#\cR$ follows. 
Let 
\[
\log : (\bR^\times)^3\ni{}^t(z_1,z_2,z_3)\mapsto {}^t(\log|z_1|,\log|z_2|,\log|z_3|)\in\bR^3
\]
be the logarithmic map. 
By Dirichlet's unit theorem, the set 
\[
\cE(F):=\{\log\bep\,|\,\bep={}^t(\ep,\ep^\sigma,\ep^{\sigma^2}),\ep\in\cO(F)^\times\}
\]
is a lattice of rank two in the plane 
$\iPi_{\log}:=\{{}^t(u_1,u_2,u_3)\in\bR^3\,|\,u_1+u_2+u_3=0\}$. 
We use the modified logarithmic map 
\[
\phi: (\bR^\times)^3\ni\bz={}^t(z_1,z_2,z_3)\mapsto
\bu={}^t(u_1,u_2,u_3)=\log(D^{-1/6}\bz)\in\bR^3.
\]
For $\bz=\bde(-x\mathbf{1}+y\bte)\in\cL\cap\cH$, the image 
$\bu=\phi(\bz)=\phi(\bde(-x\mathbf{1}+y\bte))$ 
is contained in the displaced lattice 
$\phi(\bde)+\cE(F)\subset\iPi_{\log}$, 
because $-x\mathbf{1}+y\bte\in\cE(F)$. 
Okazaki \cite[Theorem 4.2]{Oka02} showed that 
\begin{align}
(1-\sigma)\phi(\bde),\ (1-\sigma^2)\phi(\bde),\ 3\,\phi(\bde)\in\cE(F)\label{m3}.
\end{align}
Note that $3\,\phi(\bde)=(1-\sigma)\phi(\bde)+(1-\sigma^2)\phi(\bde)\in\cE(F)$. 
Hence 
\[
\cM=\bZ\,\phi(\bde)+\cE(F)\subset\iPi_{\log}
\]
is a lattice with discriminant $d(\cM)=d(\cE(F))$ or $\frac{1}{3}d(\cE(F))$.

The curve $\cH$ consists of three connected components $\cH\cap\{z_k>0\}$, $(k=1,2,3)$. 
We divide $\cH$ in another way into three parts 
\[
\cH_k=\cH\cap\{|z_k|\leq |z_{k+1}|,|z_{k+2}|\},\quad (k=1,2,3)
\]
where we take the subscripts of $z_k$s modulo $3$. 
We put 
\[
\cC:=\phi(\cH)\subset\iPi_{\log},\quad \cC_k:=\phi(\cH_k),\quad (k=1,2,3)
\]
(see also \cite[Figure 1]{Oka02}). 
For example, the curve $\cC_1$ is given as 
\[
\cC_1 : u_1+u_2+u_3=0,\quad e^{u_1}=\pm(e^{u_2}-e^{u_3}),\quad u_1\leq u_2,u_3.
\]
The map $\phi$ is a bijection from $\cH$ to $\cC$ (resp. $\cH_k$ to $\cC_k$). 
For $\bz=\bde(-x\mathbf{1}+y\bte)\in\cL\cap\cH$, we get 
$\bz\in\cL\cap\cH_k$, 
$\bz^\sigma\in\cL\cap\cH_{k+1}$, 
$\bz^{\sigma^2}\in\cL\cap\cH_{k+2}$ for some $k$, and hence  
$\bu=\phi(\bz)\in \cM\cap\cC_k$, $\phi(\bz^\sigma)\in\cM\cap\cC_{k+1}$, 
$\phi(\bz^{\sigma^2})\in\cM\cap\cC_{k+2}$ for some $k$. 

We adopt local coordinates for each piece $C_k\subset\iPi_{\log}$ by 
\begin{align}
s=s(\bu):=\frac{u_{k+1}-u_{k+2}}{\sqrt{2}},\quad 
t=t(\bu):=-\frac{\sqrt{6}u_k}{2}.\label{defst}
\end{align}

Then we get the following exponential gap principle (cf. also \cite{Wak07a}, \cite{Akh09}): 
\begin{theorem}[{\cite[Theorem 5.6]{Oka02}}]\label{thOkaG}
Let $\cM$ be a lattice of rank $2$ in $\iPi_{\rm log}$. 
Assume distinct points $\bu$ and $\bu'$ of $\cM$ lie on the same piece 
$\cC_k$ of the curve $\cC$. 
Set $t=t(\bu)$ and $t'=t(\bu')$. 
Assume $t'\geq t$. 
Then we have 
\[
t'\geq\frac{\sqrt{2}\,d(\cM)\,\mathrm{exp}(\sqrt{6}t/2)}
{1+\mathrm{exp}(-2(t'-t)/\sqrt{6}\log 2)}.
\]
Moreover, this gap principle can be rewritten in terms of $r=||\bu||$ and $r'=||\bu'||$ 
since $r'\geq t'$ and $t\geq r/2$ $($or $t\geq r/1.01$ if $r\geq 1.2)$.
\end{theorem}

Now we take the cyclic cubic field $L:=\bQ(\theta_1)$ and the full unit group 
$\mathcal{O}_L^\times$. 
Let 
\[
\cE(L):=\{\log\bep\,|\,\bep={}^t(\ep,\ep^\sigma,\ep^{\sigma^2}),\ep\in\cO_L^\times\}
\]
be a lattice on $\iPi_{\log}$. 
Okazaki also gave the following upper bound by using a result of Laurent-Mignotte-Nesterenko 
\cite{LMN95}. 
\begin{theorem}[{\cite[Theorem 7.2]{Oka02}}]\label{thOkaG2}
For $\bz'\in\cL\cap\cH$ and $t'=t(\bz')$, we have 
\[
\frac{t'}{d(\bZ\,\phi(\bde)+\cE(L))}\leq 5.04\cdot 10^4.
\]
\end{theorem}

In order to give a proof of Theorem \ref{thOkab}, we recall a proof of Theorem \ref{thOkaM} 
(see \cite[page 309]{Oka02}). 
Assume that $\#\cR>3$. 
Then there exist distinct points $\bu$, $\bu'\in\cM\cap\cC_k$ for some $k$ 
and we may apply Theorem \ref{thOkaG} to $\cM=\bZ\,\phi(\bde)+\cE(L)$. 
By \cite[Lemma 5.7]{Oka02}, we see 
\begin{align}
d(\cE(L))=\sqrt{3}R\geq 0.6\label{es}
\end{align}
where $R$ is the regulator of $L$. 
Hence we obtain $d(\cM)\geq \frac{1}{3}d(\cE(L))\geq 0.2$. 
Set $t:=t(\bu)$, $t':=t(\bu')$ with $t\leq t'$. 
By Theorem \ref{thOkaG} and $t'-t\geq 0$, we get 
$t'\geq \sqrt{2}\cdot(0.2)\cdot\exp(\sqrt{6}t/2)/(1+1)>0.14\exp(\sqrt{6}t/2)$. 
Hence 
\[
0.14\exp(\sqrt{6}t/2)-t < t'-t.
\]
By Theorems \ref{thOkaG} and \ref{thOkaG2}, we get 
\begin{align}
&\frac{\sqrt{2}\,\mathrm{exp}(\sqrt{6}t/2)}
{1+\mathrm{exp}(-2(0.14\exp(\sqrt{6}t/2)-t)/\sqrt{6}\log 2)}\label{eqq}\\
&<\frac{\sqrt{2}\,\mathrm{exp}(\sqrt{6}t/2)}
{1+\mathrm{exp}(-2(t'-t)/\sqrt{6}\log 2)}\leq 
\frac{t'}{d(\cM)}\leq 5.04\cdot 10^4\nonumber
\end{align}
and hence $t\leq 8.56$. 
We also have $r=||\bu||\leq 1.01\cdot 8.56<8.65$ (see Theorem \ref{thOkaG}). 

On the other hand, by \cite[Lemma 5.8]{Oka02}, 
if $D=D(F)>10^{12}$ then 
\[
||\bu||\geq\frac{1}{2\sqrt{6}}\log\frac{D}{1.01}.
\]
Hence it follows from 
\[
8.65>\frac{1}{2\sqrt{6}}\log\frac{D}{1.01}
\]
that $D<2.56\cdot 10^{18}$. 
This is a sketch of the proof of Theorem \ref{thOkaM}. \\

%%%%%%%%%%%%%%%%%%%%%%%%%%%%%%%%%%%%%%%%%%%%%%%%%%%%%%%%%%%%%%%%%%%%%%%%%%%%%%%%%%%%
\begin{proof}[Proof of Theorem \ref{thOkab}]
We consider the case where $F=F_m$. 
For $m\in\bZ$, we take 
\begin{align*}
F_m(X,Y)&=(X-\te_1^{(m)}Y)(X-\te_2^{(m)}Y)(X-\te_3^{(m)}Y),
\end{align*}
and $L_m=\bQ(\te_1^{(m)})$. 
We may assume that 
\begin{align}
\te_2^{(m)}=\frac{1}{-\te_1^{(m)}-1},\quad 
\te_3^{(m)}=\frac{-\te_1^{(m)}-1}{\te_1^{(m)}}\label{te23}
\end{align}
and
\begin{align}
-2<\te_3^{(m)}<-1,\quad -\frac{1}{2}<\te_2^{(m)}<0,\quad 1<\te_1^{(m)},\label{iv}
\end{align}
because $\mathrm{N}(\te_1^{(m)})=\te_1^{(m)}\te_2^{(m)}\te_3^{(m)}=1$.
We assume that $L_m=L_n$ for $-1\leq m<n$. 
For $l=m$, $n$, we take a vector 
$\bte^{(l)}={}^t(\te_1^{(l)},\te_2^{(l)},\te_3^{(l)})\in (L_l)^3$ and put 
\begin{align*}
\bde^{(l)}&:=\mathbf{1}\times\bte^{(l)}
={}^t(\te_2^{(l)}-\te_3^{(l)},\te_3^{(l)}-\te_1^{(l)},\te_1^{(l)}-\te_2^{(l)}). 
\end{align*}
Take a common trivial solution $(x,y)=(1,0)\in\cR(F_m)\cap\cR(F_n)$ to $F_m(x,y)=F_n(x,y)=1$. 
Then, for $l=m$, $n$, we have 
\begin{align*}
\bz^{(l)}&=\bde^{(l)}(-x\mathbf{1}+y\bte^{(l)})=-\bde^{(l)}
={}^t(\te_3^{(l)}-\te_2^{(l)},\te_1^{(l)}-\te_3^{(l)},\te_2^{(l)}-\te_1^{(l)})\in 
\cL_l\cap\cH_1^{(l)}
\end{align*}
where $\cL_l=\bde^{(l)}(\bZ\mathbf{1}+\bZ\bte^{(l)})$ and $\cH_1^{(l)}$ is one of 
three parts of $\cH^{(l)}$:  
\begin{align*}
\cH_1^{(l)} &: z_1+z_2+z_3=0,\ |z_1|\leq|z_2|,|z_3|,\ z_1z_2z_3=\sqrt{D(F_l)}=l^2+3l+9,
\end{align*}
because $\Bigl{|}\frac{\te_3^{(l)}-\te_2^{(l)}}{\te_1^{(l)}-\te_3^{(l)}}\Bigr{|}<1$ and 
by (\ref{iv}), 
$\Bigl{|}\frac{\te_2^{(l)}-\te_3^{(l)}}{\te_1^{(l)}-\te_2^{(l)}}\Bigr{|}
=\Bigl{|}\frac{1}{\te_1^{(l)}}\Bigr{|}<1$. 

For $l=m,n$, we take 
$\bu^{(l)}=\phi_l(\bz^{(l)})=\phi_l(\bde^{(l)})$ where 
$\phi_l(\bz^{(l)})=\log(D(F_l)^{-1/6}\bz^{(l)})$. 
Then we see $\bu^{(l)}\in \cC_1$ where 
\[
\cC_1 : u_1+u_2+u_3=0,\quad e^{u_1}=\pm(e^{u_2}-e^{u_3}),\quad u_1\leq u_2,u_3.
\]
We have $\bu^{(l)}\in\cM_l=\bZ\phi_l(\bde_l)+\cE(L_l)$. 
It follows from (\ref{m3}) that $3\phi_m(\bde^{(m)})$, $3\phi_n(\bde^{(n)})
\in\cE(L_m)=\cE(L_n)$. 
We consider the lattice 
\[
\cM_{m,n}:=\bZ\phi_m(\bde^{(m)})+\bZ\phi_n(\bde^{(n)})+\cE(L_n)
\]
of rank two with $[\cM_{m,n}:\cE(L_n)]=1$, $3$ or $9$. 
Set $t:=t(\bu^{(m)})=-\sqrt{6}u_1^{(m)}/2$ and $t':=t(\bu^{(n)})=-\sqrt{6}u_1^{(n)}/2$ 
as in (\ref{defst}). 
We see that $t< t'$. 
Now we may apply Theorem \ref{thOkaG} to $\cM=\cM_{m,n}$ and two points 
$\bu^{(m)},\bu^{(n)}\in\cM_{m,n}$. 

(i) The case where $[\cM_{m,n}:\cE(L_n)]\leq 3$. 
We get $t\leq 8.56$ and $m\leq 35731$ as in (\ref{eqq}). 

(ii) The case where $[\cM_{m,n}:\cE(L_n)]=9$. 
By the same argument as in (\ref{eqq}), it follows from 
Theorems \ref{thOkaG} and \ref{thOkaG2} that 
\begin{align}
&\frac{\sqrt{2}\,\mathrm{exp}(\sqrt{6}t/2)}
{1+\mathrm{exp}(-2(0.04\exp(\sqrt{6}t/2)-t)/\sqrt{6}\log 2)}
<\frac{t'}{d(\cM_{m,n})}\leq3\cdot5.04\cdot 10^4.\label{eqq2}
\end{align}
We get $t\leq 9.46$ and hence $m\leq 107588$. 
We consider the index $j_u:=[\cO_{L_m}^\times : 
\langle-1,\theta^{(m)}_1,\theta^{(m)}_2\rangle]$. 
If $m^2+3m+9$ is squarefree then $j_u=1$ (see \cite{Tho79}, \cite{Was87}). 
Ennola \cite{Enn91} confirmed that for $-1\leq m\leq 10000$, 
$j_u=1$ except for $m=3,5,12,54,66,1259,2389$. 
Using a computer software MAGMA \cite{BCP97}, 
we may check that $j_u=1$ also for $10000<m\leq 107588$. 

For such $m$, we have $\cE(L_m)=\cE(L_n)=\bZ\mathbf{e}_1+\bZ\mathbf{e}_2$ where 
$\mathbf{e}_1=\log \bte^{(m)}$ and $\mathbf{e}_2=\sigma \mathbf{e}_1$. 
Note that $\sigma \mathbf{e}_2=-\mathbf{e}_1-\mathbf{e}_2$ by (\ref{te23}).
We write $\phi_l(\bde^{(l)})=\frac{a}{3}\mathbf{e}_1+\frac{b}{3}\mathbf{e}_2$ 
with $a,b\in\bZ$. 
From (\ref{m3}), $(1-\sigma)\phi_l(\bde^{(l)})=\frac{a+b}{3}\mathbf{e}_1
+\frac{-a+2b}{3}\mathbf{e}_2\in\cE(L_n)$. 
Hence $(a,b)\equiv (0,0),(1,2),(2,1)\pmod{3}$. 
This implies that $[\cM_{m,n}:\cE(L_n)]\leq 3$ and it 
contradicts the assumption $[\cM_{m,n}:\cE(L_n)]=9$.
\end{proof}
\begin{remark}
We may check the condition $3\phi_m(\bde^{(m)})\in\cE(L_m)$ directly as follows: 
The elements 
$3\phi_m(\bde^{(m)})=\log\bep$, $\bep={}^t(\ep,\ep^\sigma,\ep^{\sigma^2})$ 
become roots of 
\[
g_m(X)=X^3+3X^2-(m^2+3m+6)X+1
\]
with discriminant $(2m+3)^2(m^2+3m+9)^2$. 
It follows from the transformation formula in \cite[page 97]{HM09b} and Theorem \ref{thC3} that 
$\mathrm{Spl}_\bQ g_m(X)=\mathrm{Spl}_\bQ f^{C_3}_{-3(m-3)/(2m+3)}(X)=L_m$. 
\end{remark}

\section{Proof of Theorem \ref{thallsol}}\label{seThue}

Let $\theta_2$ be a root of $f_m^{C_3}(X)$ with $-\frac{1}{2}<\theta_2<0$. 
For $m\geq 18$, we have 
\begin{align*}
-\frac{1}{m}+\frac{2}{m^2}-\frac{3}{m^3}-\frac{3}{m^4}+\frac{17}{m^5}-\frac{28}{m^6}
<\te_2<-\frac{1}{m}+\frac{2}{m^2}-\frac{3}{m^3}-\frac{3}{m^4}+\frac{17}{m^5}-\frac{27}{m^6}.
\end{align*}
Let $[\alpha]$ be the greatest integer less than or equal to $\alpha$. 
\begin{lemma}\label{lemconv}
The continued fraction expansion of the root $\te_2$ of $f_m^{C_3}(X)$ is given 
as follows:\\
When $m$ is even, 
\begin{align*}
\te_2=\begin{cases}
[-1;1,m+1,\frac{m}{2},1,3,\frac{m-14}{14},1,6,[\frac{m}{6}],\ldots]\ \hspace*{12mm}
\ \mathrm{if}\ m=14k\geq 28,\\
[-1;1,m+1,\frac{m}{2},1,3,\frac{m-2}{14},[\frac{49m+72}{6}],\ldots]\ \hspace*{13mm}
\ \mathrm{if}\ m=14k+2\geq 156,\\
[-1;1,m+1,\frac{m}{2},1,3,\frac{m-4}{14},6,1,[\frac{m-4}{6}],\ldots]\ \hspace*{9.6mm}
\ \mathrm{if}\ m=14k+4\geq 18,\\
[-1;1,m+1,\frac{m}{2},1,3,\frac{m-6}{14},3,2,[\frac{m-2}{6}],\ldots]\ \hspace*{9.6mm}
\ \mathrm{if}\ m=14k+6\geq 20,\\
[-1;1,m+1,\frac{m}{2},1,3,\frac{m-8}{14},2,3,[\frac{m-2}{6}],\ldots]\ \hspace*{9.6mm}
\ \mathrm{if}\ m=14k+8\geq 22,\\
[-1;1,m+1,\frac{m}{2},1,3,\frac{m-10}{14},1,1,2,1,[\frac{m-4}{6}],\ldots]\  
\ \mathrm{if}\ m=14k+10\geq 24,\\
[-1;1,m+1,\frac{m}{2},1,3,\frac{m-12}{14},1,2,1,1,[\frac{m-2}{6}],\ldots]\  
\ \mathrm{if}\ m=14k+12\geq 68.
\end{cases}
\end{align*}
When $m$ is odd, 
\begin{align*}
\te_2=\begin{cases}
[-1;1,m+1,\frac{m+1}{2},3,1,\frac{m-15}{14},2,3,[\frac{m-1}{6}],\ldots] \hspace*{9.3mm}
\ \mathrm{if}\ m=14k+1\geq 29,\\
[-1;1,m+1,\frac{m+1}{2},3,1,\frac{m-17}{14},1,1,2,1,[\frac{m-3}{6}],\ldots]\ 
\ \mathrm{if}\ m=14k+3\geq 31,\\
[-1;1,m+1,\frac{m+1}{2},3,1,\frac{m-19}{14},1,2,1,1,[\frac{m-3}{6}],\ldots]\ 
\ \mathrm{if}\ m=14k+5\geq 33,\\
[-1;1,m+1,\frac{m+1}{2},3,1,\frac{m-21}{14},1,6,[\frac{m-1}{6}],\ldots] \hspace*{9.6mm}
\ \mathrm{if}\ m=14k+7\geq 35,\\
[-1;1,m+1,\frac{m+1}{2},3,1,\frac{m-9}{14},2,3,[\frac{49m+73}{6}],\ldots] \hspace*{6.7mm}
\ \mathrm{if}\ m=14k+9\geq 65,\\
[-1;1,m+1,\frac{m+1}{2},3,1,\frac{m-11}{14},6,1,[\frac{m-5}{6}],\ldots] \hspace*{9.6mm}
\ \mathrm{if}\ m=14k+11\geq 25,\\
[-1;1,m+1,\frac{m+1}{2},3,1,\frac{m-13}{14},3,2,[\frac{m-3}{6}],\ldots] \hspace*{9.6mm}
\ \mathrm{if}\ m=14k+13\geq 27.
\end{cases}
\end{align*}
\end{lemma}

\begin{theorem}[{\cite[Theorem 3]{LPV99}}]\label{thLPV}
Let $m\geq 1$ and assume that $(x,y)\in\bZ^2$ is a primitive solution to 
$|F_m(x,y)|\leq \lambda(m)$ with $-\frac{y}{2}<x\leq y$ and $\frac{8\lambda(m)}{2m+3}\leq y$ 
where $\lambda(m):\bZ\rightarrow \bN$. Then\\
{\rm (i)} $x/y$ is a convergent to $\theta_2$, and we have either $y=1$ or 
\[
\left|\frac{x}{y}-\theta_2\right|<\frac{\lambda(m)}{y^3(m+1)}
\quad\mathrm{and}\quad y\geq m+2.
\]
{\rm (ii)} Put
\[
\kappa=\frac{\log(\sqrt{m^2+3m+9})+0.83}{\log(m+\frac{3}{2})-1.3}.
\]
If $m\geq 30$ we have 
\begin{align}
y^{2-\kappa}<17.78\cdot 2.59^\kappa\lambda(m).\label{eqy2}
\end{align}
\end{theorem}
%

%%%%%%%%%%%%%%%%%%%%%%%%%%%%%%%%%%%%%%%%%%%%%%%%%%%%%%%%%%%%%%%%%%%%%%%%%%%%%%%%%%%%%%%
\begin{proof}[Proof of Theorem \ref{thallsol}]

By Theorem \ref{thmain} and Theorem \ref{thOkab}, it is enough to find all non-trivial 
solutions $(x,y)\in\bZ^2$ to $F_m(x,y)=\lambda$ with $\lambda\mid m^2+3m+9$ 
for $-1\leq m\leq 35731$. 
Indeed if there exists a non-trivial solution $(x,y)\in\bZ^2$ to 
$F_n(x,y)=\lambda\mid n^2+3n+9$ for $n\geq 35732$ 
then there exists $-1\leq m\leq 35731$ such that $L_m=L_n$ by Theorem \ref{thOkab}. 
By Theorem \ref{thmain}, such $n$ can be found by solutions 
$(x,y)\in\bZ^2$ to $F_m(x,y)=\lambda$ via (\ref{eqth2}).\\ \\
(i) The case where $m=-1$. 
Peth\"o \cite{Pet87} solved the Thue inequalities 
$|F_{-1}(x,y)|\leq 200$ with $|y|\leq 10^{500}$ and de Weger \cite{Weg95} 
determined the complete solutions to $|F_{-1}(x,y)|\leq 10^6$. 
Only nine solutions to these satisfy that $\lambda>0$ divides $m^2+3m+9=7$ as on Table $1$. \\ \\
(ii) The case where $0\leq m\leq 2407$. 
For small $m$, we may use a result of Bilu-Hanrot \cite{BH96}. 
By using the following command of MAGMA \cite{BCP97} 
which is based on \cite{BH96}, we confirmed that all non-trivial solutions to (\ref{eqFth}) 
are given in Table $1$. 
\begin{verbatim} 
R<x>:=PolynomialRing(Integers());
for m:=0 to 2407 do
d:=Divisors(m^2+3*m+9);
 for i:=1 to #d do 
 s:=Solutions(Thue(x^3-m*x^2-(m+3)*x-1),d[i]); 
 nts:=[s[i] : i in [1..#s] | s[i,1]*s[i,2]*(s[i,1]+s[i,2]) ne 0];
  if #nts gt 0 then print "Non-trivial", [m,d[i]], nts; 
  end if;
 end for;
end for;
\end{verbatim}
(iii) The case where $2408\leq m\leq 35731$ and $2(2m+3+\frac{27}{2m+3})\leq y$. 
We consider the Thue inequality 
\begin{align}
|F_m(x,y)|\leq m^2+3m+9\label{thueine}. 
\end{align}
If $(x,y)\in\bZ^2$ is a solution to $(\ref{thueine})$ then 
$\pm(x,y),\pm(y,-x-y),\pm(-x-y,x)$ are also solutions to $(\ref{thueine})$. 
Hence we may assume that $-\frac{y}{2}<x\leq y$ without loss of generality. 
Applying Theorem \ref{thLPV} to 
\[
\lambda(m)=m^2+3m+9,\quad 
\frac{8\lambda(m)}{2m+3}=2\left(2m+3+\frac{27}{2m+3}\right), 
\]
$x/y$ is a convergent to $\te_2$. 
Write $m=14k+l$ with $k\in\bZ$ and $0\leq l\leq 13$. 
For $0\leq l\leq 13$, we define $i_l\in\bZ$ by the following table: 
\begin{center}
\begin{tabular}{|c||c|c|c|c|c|c|c|c|c|c|c|c|c|c|}\hline
$l$ & $0$ & $1$ & $2$ & $3$ & $4$ & $5$ & $6$ & $7$ & $8$ & $9$ & $10$ 
& $11$ & $12$ & $13$\\\hline
$i_l$ & $9$ & $9$ & $7$ & $11$ & $9$ & $11$ & $9$ & $9$ & $9$ & $7$ & $11$ 
& $9$ & $11$ & $9$\\\hline
\end{tabular}
\end{center}
Then by Lemma \ref{lemconv}, for each $0\leq l\leq 13$, 
we may show that the denominator $q_{i_l}$ of the ${i_l}$-th convergent 
of $\te_2$ is a polynomial 
of degree four with respect to $m$ and the inequality $(\ref{eqy2})$ 
does not hold for $m\geq 2408$. 
 
For example, if $l=0$ then the denominators of the first convergents of $\te_2$ is given 
as $1$, $m+2$, $\frac{1}{2}(m^2+2 m+2)$, $\frac{1}{2}(m^2+4 m+6)$, $2m^2+7m+10$, 
$\frac{1}{7}(m^3-7m^2-30m-49)$, $\frac{1}{7}(m+3)(m^2+4m+7)$, 
$m^3+5m^2+12m+11$, 
\begin{align*}
q_9=\begin{cases}
\frac{1}{42}(7m^4+41m^3+126m^2+191m+126)\ \ \mathrm{if}\ m=42k',\\
\frac{1}{42}(7m^4+27m^3+56m^2+23m-28)\hspace*{8.8mm} \mathrm{if}\ m=42k'+14,\\
\frac{1}{42}(7m^4+13m^3-14m^2-145m-182)\hspace*{4.6mm} \mathrm{if}\ m=42k'+28,
\end{cases}
\end{align*}
and the inequality 
\[
q_9^{2-\kappa}<17.78\cdot 2.59^\kappa (m^2+3m+9)
\]
does not hold for $m\geq 2408$. 
We omit to explain the case where $1\leq l\leq 13$. 

For $2408\leq m\leq 35731$, we checked that the $i$-th convergent 
$p_i/q_i$ with $2(2m+3+\frac{27}{2m+3})\leq q_i$ $(1\leq i\leq i_l-1)$
does not satisfy (\ref{eqFth}). \\ \\
(iv) The case where $2408\leq m\leq 35731$ and $y< 2(2m+3+\frac{27}{2m+3})$. 
If $(x,y)\in\bZ^2$ is a solution to $(\ref{eqFth})$ then 
$(y,-x-y)$ and $(-x-y,x)$ are also solutions to $(\ref{eqFth})$. 
Hence we may assume that $y\geq 1$. 
The bound is small enough to reach using a computer. 
Indeed we checked that for $2408\leq m\leq 35731$, $\lambda>0$ with $\lambda\mid m^2+3m+9$ 
and $1\leq y< 2(2m+3+\frac{27}{2m+3})$, every solution $x\in\bZ$ to 
$F_m(x,y)-\lambda=0$ satisfies $x=-y$. 
Thus there are no further non-trivial solutions. 
\end{proof}
\begin{proof}[Proof of Corollary \ref{corsol}]
We should determine the all trivial solutions $(x,y)\in\bZ^2$, 
i.e. $F_m(x,y)$ $=$ $m^2+3m+9=c^3$ with $xy(x+y)=0$ for some $c\in\bZ$. 
It follows from Lemma \ref{lemcond} that if $m^2+3m+9=c^3$ then the conductor 
$\mathfrak{f}_m$ of $L_m$ should be $9$. 
By Theorem \ref{thOka}, however, only the simplest cubic fields $L_m, (m\geq -1)$ 
of conductor $9$ are $L_0$, $L_3$ and $L_{54}$. This implies that $c=3$. 
This can be also verified by showing that the elliptic curve $y^2+3y=x^3-9$ over $\bQ$ 
has only the two integral points $(x,y)=(3,-6)$, $(3,3)$. 
\end{proof}
\begin{acknowledgments}
The author thanks Professor Ryotaro Okazaki for informing him Theorem \ref{thOka} 
and for explaining an idea of the proof of Theorem \ref{thOkab}. 
He also thanks Professor Isao Wakabayashi for helpful discussions. 
This work was partially supported by Rikkyo University Special Fund for Research 
and by the Grant-in-Aid for Young Scientists (B) No. 22740028, 
The Ministry of Education, Culture, Sports, Science and Technology, Japan.
\end{acknowledgments}

%%%%%%%%%%%%%%%%%%%%%%%%%%%%%%%%%%%%%%%%%%%%%%%%%%%%%%%%%%%%%%%%%%%%%%%%%

\vspace*{2mm}
\noindent
Akinari Hoshi\\
Department of Mathematics\\
Rikkyo University\\ 
3--34--1 Nishi-Ikebukuro Toshima-ku\\ 
Tokyo, 171--8501, Japan\\
E-mail: \texttt{hoshi@rikkyo.ac.jp}\\
Web: \texttt{http://www2.rikkyo.ac.jp/web/hoshi/}

\end{document}